\documentclass[reqno,centertags, 12pt]{amsart}
\usepackage{amsmath,amsthm,amscd,amssymb}
\usepackage{latexsym}
\newcommand{\bbR}{{\mathbb{R}}}

\newcommand{\bbZ}{{\mathbb{Z}}}

\newcommand{\calI}{{\mathcal I}}

\newcommand{\calT}{{\mathcal T}}

\newcommand{\lb}{\label}
\newcommand{\f}{\frac}

\newcommand{\tr}{\text{\rm{Tr}}}

\newcommand{\dist}{\text{\rm{dist}}}
\newcommand{\loc}{\text{\rm{loc}}}

\newcommand{\rank}{\text{\rm{rank}}}

\newcommand{\bi}{\bibitem}

\newcommand{\beq}{\begin{equation}}
\newcommand{\eeq}{\end{equation}}
\newcommand{\ba}{\begin{align}}
\newcommand{\ea}{\end{align}}
\newcommand{\veps}{\varepsilon}

\newcounter{smalllist}

\allowdisplaybreaks
\numberwithin{equation}{section}

\newtheorem{theorem}{Theorem}[section]

\newtheorem*{CON}{Conjecture}
\newtheorem*{p2.1}{Proposition 2.1}
\newtheorem{proposition}[theorem]{Proposition}
\newtheorem{lemma}[theorem]{Lemma}

\theoremstyle{definition}

\theoremstyle{remark}

\newcommand{\abs}[1]{\lvert#1\rvert}

\begin{document}
\title[Eigenvalue Bounds]{Eigenvalue Bounds in the Gaps of Schr\"odinger
Operators and Jacobi Matrices}
\author[D.~Hundertmark and B.~Simon]{Dirk Hundertmark$^1$ and Barry Simon$^2$}

\thanks{$^1$ School of Mathematics, Watson Hall, University of Birmingham,
Birmingham, B15 2TT, UK. On leave from Department of Mathematics, Altgeld Hall,
University of Illinois at Urbana-Champaign, Urbana, IL 61801, USA.
Email: hundertd@for.mat.bham.ac.uk. Supported in part by NSF grant
DMS-0400940. }

\thanks{$^2$ Mathematics 253-37, California Institute of Technology, Pasadena, CA 91125, USA.
E-mail: bsimon@caltech.edu. Supported in part by NSF grant DMS-0140592 and
U.S.--Israel Binational Science Foundation (BSF) Grant No.\ 2002068}

\date{May 18, 2007}

\begin{abstract} We consider $C=A+B$ where $A$ is selfadjoint with a gap $(a,b)$ in
its spectrum and $B$ is (relatively) compact. We prove a general result allowing
$B$ of indefinite sign and apply it to obtain a $(\delta V)^{d/2}$ bound for
perturbations of suitable periodic Schr\"odinger operators and a (not quite)
Lieb--Thirring bound for perturbations of algebro-geometric almost periodic Jacobi
matrices.
\end{abstract}

\maketitle

\section{Introduction} \lb{s1}

The study of the eigenvalues of Schr\"odinger operators below the essential spectrum
goes back over fifty years to Bargmann \cite{Barg3}, Birman \cite{Bir61}, and
Schwinger \cite{Schw}, and of power bounds on the eigenvalues to Lieb--Thirring
\cite{LT1,LT2}.

There has been considerably less work on eigenvalues in gaps---much of what has been
studied followed up on seminal work by Deift and Hempel \cite{DH86}; see
\cite{AADH,ADH,GGHK,GS,Hem89,Hem92,Hem97,Klaus,Lev,Saf98,Saf01JMAA,Saf01}
and especially work by Birman and collaborators
\cite{Bir90,Bir91,Bir91-ASM,Bir91-FAA,Bir95,Bir97,Bir98,BLS,BP,BR,BW}.
Following Deift--Hempel, this work has mainly focused on the set of $\lambda$'s
so that some given fixed $e$ in a gap of $\sigma(A)$ is an eigenvalue of $A+\lambda B$
and the growth of the number of eigenvalues as $\lambda\to\infty$ most often for closed
intervals  strictly inside the gap. Most, but not all, of this work has focused on $B$'s
of a definite sign. Our goal in this note is to make an elementary observation
that, as regards behavior at an edge for fixed $\lambda$, allows perturbations of
either sign. The decoupling in steps we use does not work for the question raised
by Deift--Hempel, which may be why it does not seem to be in the literature.

We will present two applications: a Cwikel--Lieb--Rozenblum-type finiteness result
\cite{Cwi,Lieb,Roz} for suitable gaps in $d\geq 3$ periodic Schr\"odinger operators
and a critical power estimate on eigenvalues in some one-dimensional almost
periodic problems.

To state our results precisely, we need some notation. For any selfadjoint operator
$C$, $E_\Omega(C)$ will denote the spectral projections for $C$. We define
\begin{equation} \lb{1.1}
\#(C\in\Omega) =\dim (E_\Omega(C))
\end{equation}
and
\begin{equation} \lb{1.2}
\#(C>\alpha) = \dim (E_{(\alpha,\infty)}(C))
\end{equation}
and similarly for $\#(C\geq\alpha)$, $\#(C<\alpha)$, $\#(C\leq\alpha)$.

We will write
\begin{equation} \lb{1.3}
B=B_+-B_-
\end{equation}
with $B_\pm\geq 0$. While often we will take $B_\pm =\max(\pm B,0)$, we do not require
$B_+B_- =0$ or $[B_+,B]=0$. Our main technical result, which we will prove in Section~\ref{s2},
is

\begin{theorem}\lb{T1.1} Let $A$ be a selfadjoint operator and $x,y\in\bbR$ so $(x,y)\cap
\sigma(A)=\emptyset$. Let $B$ be given by \eqref{1.3} with $B_+,B_-$ both compact. Let
$C=A+B$. Let $x < e_0 < e_1 =\f12 (x+y)$, then
\begin{equation} \lb{1.4}
\#(C\in (e_0, e_1))\leq \#(B_+^{1/2} (e_0-A)^{-1} B_+^{1/2}\geq 1) +
\#(B_-\geq\tfrac12 (y-x))
\end{equation}
\end{theorem}

In Section~\ref{s3}, we discuss an analog when $A$ is unbounded but bounded below and
$B_\pm$ are only relatively compact.

If $V$ is a periodic locally $L^{d/2}$ function on $\bbR^d$ ($d\geq 3$), then $A=-\Delta
+V$ can be written as a direct integral of operators, $A(k)$, with compact resolvent,
with the integral over the fundamental cell of a dual lattice (see \cite{RS4}). If
$\veps_1(k) \leq \veps_2 (k) \leq\dots$ are the eigenvalues of $A(k)$, then $(x,y)$ is
a gap in $\sigma(A)$ (i.e., connected component of $\bbR\setminus\sigma(A)$) if and only
if there is $\ell$ with
\begin{equation} \lb{1.5}
\max_k\, \veps_{\ell-1}(k) =x<y =\min_k \, \veps_\ell(k)
\end{equation}
We say $y$ is a nondegenerate gap edge if and only if
\begin{equation} \lb{1.6}
\min_k\, \veps_{\ell+1}(k) >y
\end{equation}
and $\veps_\ell(k)=y$ at a finite number of points $\{k_j\}_{j=1}^N$ in the unit cell so
that for some $C$ and all $k$ in the unit cell,
\begin{equation} \lb{1.7}
\veps_\ell (k)-y\geq C\min\abs{k-k_j}^2
\end{equation}
There is a similar definition at the bottom edge if $x>-\infty$. It is a general theorem
\cite{S181} that the bottom edge is always nondegenerate. In Section~\ref{s4}, we will prove

\begin{theorem}\lb{T1.2} Let $d\geq 3$. Let $V\in L_\loc^{d/2}(\bbR^d)$ be periodic and let
$W\in L^{d/2}(\bbR^d)$. Let $(x,y)$ be a gap in the spectrum $A=-\Delta+V$ which is
nondegenerate at both ends, and let $N_{(x,y)}(W)=\#(-\Delta+V+W\in (x,y))$. Then
$N_{(x,y)}(W) <\infty$.
\end{theorem}

This will be a simple extension of the result of Birman \cite{Bir95} who proved this
if $W$ has a fixed sign. Note we have not stated a bound by $\|W\|_{d/2}^{d/2}$. This
is discussed further in Section~\ref{s4}.

In the final section, Section~\ref{s5}, we will consider
certain two-sided Jacobi matrices, $J$, on $\ell^2(\bbZ)$ with
\begin{equation} \lb{1.9}
J_{k\ell} =\begin{cases}
b_k & k=\ell \\
a_k & \ell=k+1 \\
a_{k-1} & \ell =k-1 \\
0 & \abs{\ell-k}\geq 2
\end{cases}
\end{equation}
If $E=\cup_{j=1}^{\ell+1} E_j$ is a finite union of bounded closed disjoint intervals, there is
an isospectral torus $\calT_E$ associated to $E$ of almost periodic $J$'s with $\sigma(J)= E$
(see \cite{AK,Apt,CSZ,Cr,DMN,PY,SY,Widom}). We conjecture the following:

\begin{CON} Let $J_0$ lie in some $\calT_E$. Let $J=J_0+\delta J$ be a Jacobi matrix
for which $\delta J$ is trace class, that is,
\begin{equation} \lb{1.10}
\sum_n\, \abs{\delta a_n} + \abs{\delta b_n}<\infty
\end{equation}
Then
\begin{equation} \lb{1.11}
\sum_{\lambda\in\sigma(J)\setminus E} \dist(\lambda, E)^{1/2} <\infty
\end{equation}
\end{CON}

For $e=[-2,2]$ so $J_0$ is the free Jacobi matrix with $a_n\equiv 1$, $b_n\equiv 0$,
this is a result of Hundertmark--Simon \cite{HunS}. It has recently been proven \cite{DKS2007}
for the case where $J_0$ is periodic, and it has recently been proven \cite{CLTB}
that \eqref{1.11} holds for the sum over $\lambda$'s above the top of the spectrum or
below the bottom. In Section~\ref{s5}, we will prove

\begin{theorem}\lb{T1.3} If \eqref{1.10} holds, then \eqref{1.11} holds if $\f12$
is replaced by any $\alpha >\f12$.
\end{theorem}

\begin{theorem}\lb{T1.4} If
\begin{equation} \lb{1.12}
\sum_n [\log(\abs{n}+1)]^{1+\veps} [\abs{\delta a_n} + \abs{\delta b_n}] <\infty
\end{equation}
for some $\veps >0$, then \eqref{1.11} holds.
\end{theorem}

Both the conjecture and Theorem~\ref{T1.4} are interesting because they imply that the
spectral measure obeys a Szeg\H{o} condition. This is discussed in \cite{CSZ}.

\section{Abstract Bounds in Gaps (Compact Case)} \lb{s2}

Our goal here is to prove Theorem~\ref{T1.1}. We begin by recalling the version of the
Birman--Schwinger principle for points in gaps, which is essentially the key to
\cite{AADH,ADH,Bir90,Bir91,Bir91-ASM,Bir91-FAA,Bir95,Bir97,Bir98,BLS,BP,BR,BW,
DH86,GGHK,GS,Hem89,Hem92,Hem97,Klaus,Lev,Saf98,Saf01JMAA,Saf01}:

\begin{proposition}\lb{P2.1} Let $A$ be a bounded selfadjoint operator with $(x,y)
\cap\sigma(A)=\emptyset$. Let $B$ be compact with $B\geq 0$. Let $e\in (x,y)$. Then
\begin{equation} \lb{2.1}
e\in \sigma(A+\mu B) \Leftrightarrow \mu^{-1}\in\sigma (B^{1/2} (e-A)^{-1} B^{1/2})
\end{equation}
with equal multiplicity. In particular,
\begin{equation} \lb{2.2}
\#(A+B\in (e,y))\leq \#(B^{1/2} (e-A)^{-1} B^{1/2}\geq 1)
\end{equation}
\end{proposition}

\begin{proof} This is so elementary that we sketch the proof. If for $\varphi\neq 0$,
\begin{equation} \lb{2.3}
(A+\mu B)\varphi =e\varphi
\end{equation}
then
\begin{equation} \lb{2.4}
B\varphi\neq 0
\end{equation}
since $e\notin\sigma(A)$. Moreover,
\begin{equation} \lb{2.5}
(e-A)^{-1} B\varphi =\mu^{-1}\varphi
\end{equation}
and \eqref{2.5} implies \eqref{2.3}. Thus
\begin{equation} \lb{2.6}
e\in\sigma(A+\mu B)\Leftrightarrow \mu^{-1}\in\sigma((e-A)^{-1}B)
\end{equation}
and \eqref{2.1} follows by $\sigma (CD)\setminus\{0\}=\sigma(DC)\setminus\{0\}$ (see, e.g.,
Deift \cite{Deift78}).

Since $\sigma(A+\mu B)\subset \sigma(A) + [-\mu\|B\|, \mu\|B\|]$ and discrete eigenvalues
are continuous in $\mu$ and strictly monotone by \eqref{2.4} and (see \cite{RS4})
\begin{equation} \lb{2.7}
\f{de(\mu)}{d\mu} = \langle\varphi,B\varphi\rangle
\end{equation}
eigenvalues of $A+B$ in $(x,y)$ must pass through $e$ as $\mu$ goes from $0$ to $1$ and \eqref{2.2}
follows from \eqref{2.1}. We only have inequality in \eqref{2.2} since eigenvalues can get
reabsorbed at $y$.
\end{proof}

\begin{proof}[Proof of Theorem~\ref{T1.1}] Let $C_+ =A+B_+$ so $C=C_+-B_-$.
By Proposition~\ref{P2.1}, if
\begin{equation} \lb{2.8}
n_1=\#(C_+\in (e_0,e_1)) \qquad
n_2 =\#(C_+\in (e_1,y))
\end{equation}
then
\begin{equation} \lb{2.9}
n_1 + n_2 \leq \#(B_+^{1/2} (e_0-A)^{-1} B_+^{1/2} \geq 1)
\end{equation}

By a limiting argument, we can suppose that $e_1$ is not an eigenvalue of $C_+$. Since
eigenvalues of $C_+-\mu B_-$ are strictly monotone decreasing in $\mu$, the number of eigenvalues
of $C$ in $(e_0,e_1)$ can only increase by passing through $e_1$. By repeating the argument
in Proposition~\ref{P2.1},
\begin{equation} \lb{2.10}
\#(C\in (e_0,e_1)) \leq n_1 + \#(B_-^{1/2} (C_+ -e_1)^{-1} B_-^{1/2}\geq 1)
\end{equation}

Now write
\begin{equation} \lb{2.11}
B_-^{1/2} (C_+-e_1)^{-1} B_-^{1/2} = D_1 + D_2 + D_3
\end{equation}
where $D_1$ has $E_{(-\infty, e_1)}(C_+)$ inserted in the middle, $D_2$ an $E_{(e_1,y)}(C_+)$,
and $D_3$ an $E_{[y,\infty)}(C_+)$. Since $D_1\leq 0$ and $\rank (D_2)\leq n_2$, we see
\begin{equation} \lb{2.12}
\#(B_-^{1/2} (C_+-e_1)^{-1} B_-^{1/2} \geq 1)\leq n_2 + \#(D_3\geq 1)
\end{equation}

Since $(C_+-e_1)^{-1} E_{[y,\infty)}(C_+)\leq (y-e_1)^{-1}=[\frac12(y-x)]^{-1}$, we have
\begin{equation} \lb{2.12a}
D_3 \leq [\tfrac12\, (y-x)]^{-1} B_-
\end{equation}
and thus
\begin{align}
\#(D_3\geq 1) & \leq \# ([\tfrac12\, (y-x)]^{-1} B_-\geq 1) \notag \\
&= \#(B_- \geq \tfrac12\, (y-x)) \lb{2.13}
\end{align}
\eqref{2.9}, \eqref{2.10}, \eqref{2.12}, and \eqref{2.13} imply \eqref{1.4}.
\end{proof}

\section{Abstract Bounds in Gaps (Relatively Compact Case)} \lb{s3}

In this section, we suppose $A$ is a semibounded selfadjoint operator with
\begin{equation} \lb{3.1}
q=\inf \sigma(A)
\end{equation}
We will suppose $B$ is a form-compact perturbation, which is a difference of two positive
form-compact perturbations. We abuse notation and write compact operators
\begin{equation} \lb{3.2}
B_\pm^{1/2} (A-e)^{-1} B_\pm^{1/2}
\end{equation}
for $e\notin\sigma(A)$ even though $B_\pm$ need not be operators --- \eqref{3.2}
can be defined via forms in a standard way.

In the bounded case, we only considered intervals in the lower half of a gap since
$A\to -A$, $B\to -B$ flips half-intervals. But, as has been noted in the unbounded case
(see, e.g., \cite{Bir95,Saf98}), there is now an asymmetry, so we will state separate results.
We start with the bottom half case:

\begin{theorem}\lb{T3.1} Let $A$ be a semibounded selfadjoint operator and $x,y\in\bbR$ so
$(x,y)\cap \sigma(A)=\emptyset$. Let $B=B_+-B_-$ with $B_+$ form-compact positive perturbations
of $A$. Let $C=A+B$ and $x < e_0 < e_1 =\f12 (x+y)$. Then
\begin{equation} \lb{3.3}
\begin{split}
\#(C\in (e_0,e_1)) &\leq \#(B_+^{1/2} (e_0-A)^{-1} B_+^{1/2}\geq 1) \\
&\quad + \# \biggl( B_-^{1/2} (A-q+1)^{-1} B_-^{1/2} \geq \tfrac12\, \biggl[ \f{y-x}{y-q+1}\biggr]\biggr)
\end{split}
\end{equation}
\end{theorem}

\begin{proof} We follow the proof of Theorem~\ref{T1.1} without change until \eqref{2.12a}
noting that instead
\begin{align}
(C_+-e_1)^{-1} E_{[y,\infty)}(C_+)
&\leq \f{y-q+1}{y-e_1}\, (C_+ + q+1)^{-1} \lb{3.4} \\
&\leq \f{y-q+1}{y-e_1}\, (A-q+1)^{-1} \lb{3.5}
\end{align}
since $q\leq A\leq C_+$ and
\[
\sup_{x\geq y}\, \f{x-q+1}{x-e_1}
\]
is taken at $x=y$ since $q-1 <e_1$. By \eqref{3.5},
\[
\#(D_3\geq 1) \leq \# \biggl(B_-^{1/2} (A-q+1)^{-1} B_-^{1/2} \geq \f{y-e_1}{y-q+1} \biggr)
\qedhere
\]
\end{proof}

\begin{theorem}\lb{T3.2} Let $A$ be a semibounded selfadjoint operator and $(x,y)\in\bbR$ so
$(x,y)\cap\sigma(A)=\emptyset$. Let $B=B_+-B_-$ with $B_\pm$ form-compact positive perturbations
of $A$. Let $C=A+B$ and $e_1=\f12 (x+y) <e_0 <y$. Then
\begin{equation} \lb{3.6}
\begin{split}
\#(C\in (e_1,e_0)) &\leq \#(B_-^{1/2} (A-e_0)^{-1} B_-^{1/2} \geq 1) \\
&\quad + \#(B_+^{1/2} (A-B_- - e_1)^{-1} E_{(-\infty,x)} (A-B_-) B_+^{1/2} \geq 1)
\end{split}
\end{equation}
\end{theorem}

\begin{proof} Identical to the proof of Theorem~\ref{T1.1} through \eqref{2.12a}.
\end{proof}

The second term in \eqref{3.6} is easily seen to be finite since the operator is compact.
However, any bound depends on both $B_+$ and $B_-$.

\section{$L^{n/2}$ Bounds in Gaps for Periodic Schr\"odinger Operators} \lb{s4}

Birman \cite{Bir95} proved for $V$\!, as in Theorem~\ref{T1.2}, and any $W$ that uniformly
in any gap $(x,y)$, $\sup_{\lambda\in (x,y)} \|\abs{W}^{1/2} (-\Delta +V-\lambda)^{-1}
\abs{W}^{1/2}\|_{\calI_{d/2}^w} \leq c\|W\|_{d/2}$ where $\|\cdot\|_{\calI_{d/2}^w}$ is a weak
$\calI_d$ trace class norm \cite{S73}. To be precise, in his Proposition~3.1, he
proved $\|\abs{W}^{1/2} (-\Delta +V -\lambda_0)^{-1} \abs{W}^{1/2}\|_{\calI_{d/2}}$ is finite
away from $x$ and $y$, and then in (3.15), he proved the weak estimate at the end
points. He used this to prove for $W$ of a definite sign
\begin{equation} \lb{4.1}
N_{(x,y)}(W) \leq c\int_{\bbR^d} \abs{W(z)}^{d/2}\, dz
\end{equation}
It implies relative compactness, and given Theorems~\ref{T3.1} and \ref{T3.2}, proves
Theorem~\ref{1.2}.

Note that, by Theorem~\ref{T3.1}, we get for any $x' >x$,
\begin{equation} \lb{4.2}
N_{(x',y)}(W) \leq c_{x'} \int_{\bbR^d} \abs{W(z)}^{d/2}\, dz
\end{equation}
but we do not get such a bound for $x'=x$ since there is a $W_-,W_+$ cross term in
\eqref{3.6}.

\section{Gaps for Perturbations of Finite Gap Almost
Periodic Jacobi Matrices} \lb{s5}

Our goal here is to prove Theorems~\ref{T1.3} and \ref{T1.4}.
Let
\begin{equation} \lb{5.1}
G_0 (n,m;\lambda) =\langle \delta_n, (J_0-\lambda)^{-1} \delta_m\rangle
\end{equation}
and let $(\lambda_0,\lambda_1)$ be a gap in $\sigma (J_0)$. As input, we need
two estimates for $G_0$ proven in \cite{CSZ}. First we have
\begin{equation} \lb{5.2}
\abs{G_0 (n,m;\lambda)} \leq C\dist (\lambda, \sigma(J_0))^{-1/2}
\end{equation}
uniformly in real $\lambda\notin\sigma (J_0)$ and $n$ and $m$.

To describe the other estimate, we need some notions. At a band edge, $\lambda_0$
(here and below, we study $\lambda_0$ but there is also an analysis at $\lambda_1$),
there is a unique almost periodic sequence $\{u_n (\lambda_0)\}_{n=-\infty}^\infty$
solving $(J_0-\lambda_0) u_n =0$. If $u_n=0$, we say $n$ is a resonance point. If
$u_n\neq 0$, we have a nonresonance. Since $u_n=0\Rightarrow u_{n\pm 1}\neq 0$,
we have lots of nonresonance points. Without loss, we will suppose henceforth
that $0$ is a nonresonance point. At a nonresonance point, $\lim_{\lambda\downarrow
\lambda_0} \dist (\lambda,\lambda_0)^{1/2} G_0 (n,n;\lambda)\neq 0$.

The Dirichlet Green's function is defined by
\begin{equation} \lb{5.3}
G_0^D (n,m;\lambda) = G_0 (n,m;\lambda) - G_0 (0,0;\lambda)^{-1}
G_0 (n,0;\lambda) G_0 (0,m;\lambda)
\end{equation}
Then \cite{CSZ} proves that if $0$ is a nonresonance at $\lambda_0$, then for
some small $\veps$,
\begin{align}
\lambda\in (\lambda_0, \lambda_0 +\veps)
&\Rightarrow \abs{G_0^D (n,n;\lambda)} \leq Cn \lb{5.4} \\
&\Rightarrow \abs{G_0^D (n,n;\lambda)} \leq C\abs{\lambda-\lambda_0}^{-1/2} \lb{5.5}
\end{align}

Following \cite{HunS}, we use (with $c_\pm=\max (\pm c, 0)$) with $a>0$,
\begin{equation} \lb{5.6}
\begin{pmatrix}
b&a \\ a & b
\end{pmatrix}
=\begin{pmatrix} b_+ + a & 0 \\ 0 & b_+ +a \end{pmatrix}
-\begin{pmatrix} a+b_- & -a \\ -a & a+b_- \end{pmatrix}
\end{equation}
to define $\delta J=\delta J_+ -\delta J_-$ where $\delta J_+$ is diagonal and given by
\begin{align}
(\delta J_+)_{n\,n} &= (\delta b_n)_+ + \delta a_{n-1} + \delta a_n \lb{5.7} \\
\intertext{and $(\delta J_-)$ is tridiagonal with}
(\delta J_-)_{n\,n+1} &= \delta a_n \lb{5.8a} \\
(\delta J_-)_{n\, n-1} &= \delta a_{n-1} \lb{5.8b} \\
(\delta J_-)_{n\,n} &= (\delta b_n)_- + \delta a_{n-1} + \delta a_n \lb{5.8c}
\end{align}

We also use the fact obtained via an integration by parts that if $f(\lambda_0)=0$,
$f$ continuous on $[\lambda_0, \lambda_0+\veps)$, and $C^1 (\lambda_0, \lambda_0 +\veps)$
with $f' >0$, then
\begin{equation} \lb{5.9}
\sum_{\substack{\lambda\in (\lambda_0, \lambda_0+\veps) \\
\lambda\in \sigma(J)}} f(\lambda) = \int_{\lambda_0}^{\lambda_0 +\veps}
f'(\lambda) \#(J\in (\lambda,\lambda_0 +\veps))\, d\lambda
\end{equation}

Since $f'\in L^1(\lambda_0,\lambda_0+\veps)$ and $\delta J_-$ is compact, Theorem~\ref{T1.1}
implies
\begin{equation} \lb{5.10}
\sum_{\substack{\lambda\in (\lambda_0, \lambda_0+\veps) \\
\lambda\in \sigma(J)}} f(\lambda) <\infty \Leftarrow \int_{\lambda_0}^{\lambda_0+\veps}
\# ((\delta J_+)^{1/2} (\lambda -J_0)^{-1} (\delta J_+)^{1/2} \geq 1) f'(\lambda)\,
d\lambda <\infty
\end{equation}
This leads to

\begin{proposition}\lb{P5.1} If $\delta J_\pm$ are trace class and
\begin{equation} \lb{5.11}
\int_{\lambda_0}^{\lambda_0 +\veps} f'(\lambda) \abs{\tr ((\delta J_+)^{1/2}
G_0^D (\cdot,\cdot \,;\lambda)(\delta J_+)^{1/2})}\, d\lambda <\infty
\end{equation}
then
\begin{equation} \lb{5.12}
\sum_{\substack{\lambda\in (\lambda_0, \lambda_0+\veps) \\
\lambda\in \sigma(J)}} f(\lambda) <\infty
\end{equation}
\end{proposition}

\begin{proof} $G_0-G_0^D$ is rank one and $\#(C\geq 1) \leq \|C\|_1$, so
\[
\#((\delta J_+)^{1/2} G_0 (\cdot,\cdot\, ;\lambda)(\delta J_+)^{1/2} \geq 1)
\leq 1+ \|(\delta J_+)^{1/2} G_0^D (\cdot,\cdot\, ;\lambda)(\delta J_+)^{1/2}\|_1
\]
The negative part of $G_0^D (\cdot,\cdot\, ;\lambda)$ is uniformly bounded in norm
by $\abs{a-\lambda}^{-1}$ where $a$ is either $\lambda_1$ or the unique eigenvalue
of the Dirichlet $J_0$ in $(\lambda_0-\lambda_1)$ and
\begin{align*}
\|C\|_1 &\leq \tr(C_+) + \tr (C_-) \\
&\leq \tr (C) + 2\tr (C_-)
\end{align*}
Thus \eqref{5.12} is implied by \eqref{5.10} so long as \eqref{5.11} holds.
\end{proof}

\begin{proof}[Proof of Theorem~\ref{T1.3}] By \eqref{5.5} and $\delta J_+\in\calI_1$,
we have
\[
\abs{\tr((\delta J_+)^{1/2} G_0^D (\cdot,\cdot\, ;\lambda)(\delta J_+)^{1/2})}
\leq C\abs{\lambda -\lambda_0}^{-1/2}
\]
so the integral in \eqref{5.11} is bounded by
\[
C \int_{\lambda_0}^{\lambda_0+\veps} \abs{\lambda -\lambda_0}^{\alpha -1}
\abs{\lambda -\lambda_0}^{-1/2} \, d\lambda <\infty
\]
so long as $\alpha -\f12 >0$.
\end{proof}

\begin{lemma}\lb{L5.2} For any $\alpha >0$, there is a $C$ so for all $x,y >1$,
\begin{equation} \lb{5.13}
\min (x,y) \leq C [\log (x+1)]^\alpha \, \f{y}{[\log (y+1)]^\alpha}
\end{equation}
\end{lemma}

\begin{proof} Pick $d\geq 1$ (e.g., $d=e^\alpha$), so $[\log (x+d)]^\alpha x^{-1}$ is
monotone decreasing on $[1,\infty)$. Then
\begin{equation} \lb{5.14}
\min (x,y) \leq [(\log (x+d))]^\alpha \, \f{y}{[\log (y+d)]^\alpha}
\end{equation}
If $y\leq x$, the right-hand side is bigger than $y$ and so $\min (x,y)$. If $y\geq x$,
the monotonicity shows
\[
\text{RHS} \geq [\log (x+d)]^\alpha \, \f{x}{[\log (x+d)]^\alpha} =x
\]
\eqref{5.13} follows since on $[1,\infty)$, $\f{\log (x+d)}{\log (x+1)}$ is
bounded above and below.
\end{proof}

\begin{proof}[Proof of Theorem~\ref{T1.4}] By \eqref{5.4}, \eqref{5.5}, and \eqref{5.13},
\[
\abs{G_0^D (n,n;\lambda)} \leq C\, \f{[\log (1+\abs{n})]^\alpha}{\abs{\lambda - \lambda_0}^{1/2}}
\, [\log (\lambda -\lambda_0)^{-1/2}]^{-\alpha}
\]
By \eqref{1.12}, we see
\[
\abs{\tr [(\delta J)^{1/2} G_0^D (\delta J)^{1/2}]} \leq
\f{C[\log (\lambda-\lambda_0)^{1/2}]^{-(1 +\veps)}}{(\lambda-\lambda_0)^{1/2}}
\]
Since
\[
\int_{\lambda_0}^{\lambda_0 +\veps} (\lambda-\lambda_0)^{-1}
[\log (\lambda -\lambda_0)^{-1/2}]^{-(1+\veps)}\, d\lambda <\infty
\]
the result follows.
\end{proof}

\bigskip

\end{document}